                                   %nouveau document

                             \def\flech{\rightarrow}

   \NoBlackBoxes \documentstyle{amsppt}

   \magnification=1100
    \loadbold
\topmatter
 \title  Brownian motion in Riemannian admissible complexes.
          \endtitle

\leftheadtext{Taoufik Bouziane}

\rightheadtext{Brownian motion}

\redefine\headmark#1{}
                          \author Taoufik Bouziane\endauthor

     \address Dr. Taoufik BOUZIANE, The Abdus Salam ICTP, Mathematics Section, strada costiera 11, 34014
     Trieste, Italy \endaddress

\email tbouzian\@ictp.trieste.it \endemail

\abstract

The purpose of this work is to construct a {\it Brownian motion}
with values in simplicial complexes with piecewise differential
structure. In order to state and prove the existence of such
Brownian motion, we define a family of continuous Markov processes
with values in an admissible complex; we call every process of
this family, {\it isotropic transport process}. We show that the
family of the isotropic processes contains a subsequence, which
converges weakly to a measure; we name it the {\it Wiener
measure}. Then, using the finite dimensional distributions of the
obtained Wiener measure, we construct a new admissible complex
valued continuous Markov process: the Brownian motion. We finished
with a geometric analysis of this Brownian motion, to determine
 the recurrent or transient behavior of such process.
  \endabstract


%\date March 24th 2003 \enddate

%\subjclassyear{2000}


   \endtopmatter

%   \NoBlackBoxes \documentstyle{amsppt}\hoffset 0.75 cm

               \document

     \head{ 0. Introduction.}\endhead

It has been proved in [19], [22] and [27] that, on a wide class of
Riemannian manifolds, the {\it Brownian motion} can be
approximated in law by Markov process, which generalizes the
isotropic scattering transport process on euclidean space [31]. On
the other hand, the Brownian motion was introduced as a tool to
achieve important results in Riemannian geometry and potential
theory, which is not surprising, since Brownian motion is
intimately connected with harmonic functions [15], Laplacian, and
other fundamental objects in mathematics. For instance, a complete
Riemannian manifold is hyperbolic exactly when the Brownian motion
is transient.

The purpose of this work is to consider the problem of defining
the concept of continuous random walk in the admissible Riemannian
complexes, in particular to construct a Brownian motion in
singular spaces in spite of the absence of second order
differential calculus.

The first section deals with some preliminaries on Riemannian
admissible complexes [3] [7] and finishes with a brief survey on
the theory of general Markov processes [14].

The second section is devoted firstly, to the construction of a
Markov process with values in admissible complex that we name {\it
the isotropic transport process} and secondly to show that this
latter process is a strong Markov one.

In the third section, we construct a family of isotropic process
and we show that this family contains a subsequence which
converges weakly to a measure, we call it the {\it Wiener}
measure.

The aim of the forth section is to achieve the construction of a
Brownian motion with values in an admissible complex by using the
finite dimensional distributions of the obtained Wiener measure.

Finally, in the fifth section, we close the paper  by studying the
transience/recurrence properties of the Brownian motion. In
particular, we show that, in the $2$-dimensional case if the
complex is complete simple connected of non positive curvature and
the number of branching faces is always greater or equal to three
then the Brownian motion is transient although (surprising) the
Euclidean Brownian motion in dimension $2$ is recurrent.

To our knowledge, there is an interesting study (I think it is the
unique) of M. Brin and Y. Kifer [8] on the Brownian motion in
singular spaces. In this study they consider the case of
$2$-dimensional simplicial complexes whose simplices are flat
Euclidean where they describe the Brownian motion in such complex
as the planer Brownian inside faces and, after hitting an edge,
goes into each adjacent face "with equal probability". Thus
actually, our work is the first one where it is shown the
existence of Brownian motion, and not only in the case of
$2$-dimensional complexes with flat simplices but, in the general
case of the admissible Riemannian complexes.

We notice that, the steps used for the construction of the
admissible complex-valued Brownian motion, can be extended to the
the general case of Hadamard spaces if we assume a given uniform
probability (sub-probability) measure on the link of each point of
the space.

\subhead{Acknowledgement}\endsubhead

 The author would like to express his thanks to Professor A.
 Verjovsky for encouraging him to investigate this subject
 and for his constant support over the years.

                 \head { 1.  Preliminaries. }\endhead

\subhead { 1.1 General theory [1] [2] [11] [17] [18]}\endsubhead

Let $X$ be a metric space with a metric $d$. A curve $c:I\flech X$
is called a {\it geodesic} if there is $v\geq 0$, called the
speed, such that every $t\in I$ has neighborhood $U\subset I$ with
$d(c(t_1),c(t_2))=v|t_1-t_2|$ for all $t_1,t_2\in U$. If the above
equality holds for all $t_1,t_2\in I$, then $c$ is called {\it
minimal geodesic}.

The space $X$ is called a {\it geodesic space} if every two points
in $X$ are connected by minimal geodesic. We assume from now on
that $X$ is complete geodesic space.

A triangle $\Delta$ in $X$ is a triple $(\sigma_1, \sigma_2,
\sigma_3)$ of geodesic segments whose end points match in the
usual way. Denote by $H_k$ the simply connected complete surface
of constant Gauss curvature $k$. A {\it comparison triangle}
$\bar{\Delta}$ for a triangle $\Delta \subset X$ is a triangle in
$H_k$ with the same lengths of sides as $\Delta$. A comparison
triangle in $H_k$ exists  and is unique up to congruence if the
lengths of sides of $\Delta$ satisfy the triangle inequality and,
in the case $k>0$, if the perimeter of $\Delta$ is
$<\frac{2\pi}{\sqrt k }$. Let $\bar{\Delta}= (\bar\sigma_1,
\bar\sigma_2, \bar\sigma_3)$ be a comparison triangle for
$\Delta=(\sigma_1, \sigma_2, \sigma_3)$, then for every point
$x\in \sigma_i$, $i=1,2,3$, we denote by $\bar x$ the unique point
on $\bar \sigma_i$ which lies at the same distances to the ends as
$x$.

Let $d$ denote the distance functions in both $X$ and $H_k$. A
triangle $\Delta$ in $X$ is $CAT_k$ {\it triangle} if the sides
satisfy the triangle inequality, the perimeter of $\Delta$ is
$<\frac{2\pi}{\sqrt k }$ for $k>0$, and if $d(x,y)\le d(\bar
x,\bar y)$, for every two points $x,y\in X$.

We say that $X$ has curvature at most $k$ and write $k_X\le k$ if
every point $x\in X$ has a neighborhood $U$ such that any triangle
in $X$ with vertices in $U$ and minimizing sides is $CAT_k$. Note
that we do not define $k_X$. If $X$ is Riemannian manifold, then
$k_X\le k$ iff $k$ is an upper bound for the sectional curvature
of $X$.

A geodesic space $X$ is called geodesically complete iff every
geodesic can be stretched in the two direction.

We say that a geodesic space $X$ is without conjugate points if
every two points in $X$ are connected by unique geodesic.

    \subhead { 1.2 Riemannian admissible complexes [26] [30]}\endsubhead

 Let $K$ be a locally finite simplicial complex, endowed with a
 piecewise smooth Riemannian metric $g$; i.e. $g$ is a family of
 smooth Riemannian metrics $g_\Delta$ on simplices $\Delta$ of $K$
 such that the restriction $g_\Delta|\Delta'=g_{\Delta'}$ for any
 simplices $\Delta'$ and $\Delta$ with $\Delta'\subset \Delta$.

 Let $K$ be a finite dimensional simplicial complex which is connected
locally finite. A map $f$ from $[a,b]$ to $K$ is called a broken
geodesic if there is a subdivision $a=t_0<t_1<...<t_{p+1}=b$ such
that $f([t_i,t_{i+1}])$ is contained in some cell and the
restriction of f to $[t_i,t_{i+1}]$ is a geodesic inside that
cell. Then define the length of the broken geodesic map $f$ to be
:
 $$ L(f)=\sum_{i=0}^{i=p} d(f(t_i),f(t_{i+1})).$$
 The length inside a cell is measured with respect to the metric
 of the cells.

 Then define $\tilde d(x,y)$, for every two points $x,y$ in $K$,
 to be the lower bound of the lengths of broken geodesics from $x$
 to $y$. $\tilde d$ is a pseudo-distance.

If $K$ is connected locally finite, then the space $(K,\tilde d)$
is a length space and hence a geodesic space if is complete (see
also [6])

An $l$-simplex in $K$ is called a {\it boundary simplex} if it is
adjacent to exactly one $l+1$ simplex. The complex $K$ is called
{\it boundaryless} if there are no boundary simplices in $K$.

\hskip 0 cm

 We say that the complex $K$ is {\it admissible}, if for
every connected open subset $U$ of $K$, the open set $U\setminus
\{ U\cap \{\text{the $(k-2)$}-\text{skeleton}\} \}$ is connected
($k$ is the dimension of $K$).

\hskip 0 cm

Let $x\in K$ be a vertex of $K$ so that $x$ is in the $l$-simplex
$\Delta_{l}$. We view $\Delta_{l}$ as an affine simplex in $\Bbb
R^l$, that is $\Delta _l =\bigcap_{i=0}^l H_i$, where
$H_0,H_1,...,H_l$ are closed half spaces in general position, and
we suppose that $x$ is in the topological interior of $H_0$. The
Riemannian metric $g_{\Delta_l}$ is the restriction to $\Delta_l$
of a smooth Riemannian metric defined in an open neighborhood $V$
of $\Delta_l$ in $\Bbb R^l$. The intersection
$T_x\Delta_l=\bigcap_{i=1}^l H_i \subset T_xV$ is a cone with apex
$0\in T_xV$, and $g_{\Delta_l}(x)$ turns it into an euclidean
cone. Let $\Delta_m\subset \Delta_l$ ($m<l$) be another simplex
adjacent to $x$. Then, the face of $T_x\Delta_l$ corresponding to
$\Delta_m$ is isomorphic to $T_x\Delta_m$ and we view
$T_x\Delta_m$ as a subset of $T_x\Delta_l$.

Set $T_xK =\bigcup_{\Delta_i\ni x} T_x\Delta_i$, we call it the
{\it tangent cone} of $K$ at $x$. Let $S_x\Delta_l$ denote the
subset of all unit vectors in $T_x\Delta_l$ and set $S_x=S_xK
=\bigcup_{\Delta_i\ni x} S_x\Delta_i$. The set $S_x$ is called the
{\it link} of $x$ in $K$. If $\Delta_l$ is a simplex adjacent to
$x$, then $g_{\Delta_l}(x)$ defines a Riemannian metric on the
$(l-1)$-simplex $S_x\Delta_l$. The family $g_x$ of Riemannian
metrics $g_{\Delta_l}(x)$ turns $S_x\Delta_l$ into a simplicial
complex with a piecewise smooth Riemannian metric such that the
simplices are spherical.

As it was characterized in [4], a two dimensional complete locally
finite simplicial complexes $(K,g)$ is curvature bounded by $k$
($k_K\le k$) iff the following three conditions hold:
  \roster
  \item the Gauss curvature of the open faces is bounded from above by
  $k$;
  \item for every edge $e$ of $K$, every two faces $f_1$, $f_2$ adjacent to $e$
  and every interior point $x\in e$ the sum of the geodesic
  curvatures $k_1(x)$, $k_2(x)$ of $e$ with respect to $f_1$,
  $f_2$ is nonpositive;
  \item for every vertex $x$ of $K$, every simple loop in $S_xK$ has length at least
  $2\pi$ (i.e. $S_xK$ is $CAT_1$ space).
  \endroster

\subhead{ 1.3 Liouville measure for the geodesic flow}\endsubhead

We assume that $K$ is an admissible $n$-dimensional Riemannian
complex. We denote by $K^{(i)}$ the $i$-skeleton of $K$ and $K'$
the set of points $x\in K$ such that $x$ is contained in the
interior of an $(n-1)$-simplex.

Let $x\in K'$. Then $x$ is contained in the interior of an
$(n-1)$-simplex $\Delta'$. For any $n$-simplex $\Delta$ whose
boundary $\partial\Delta$ contains $x$, let $S'_x\Delta$ denote
the open hemisphere of unit tangent vectors at $x$ pointing inside
$\Delta$. Let $\Delta_1,...,\Delta_m$, $m\ge 2$, be the
$n$-simplices containing $\Delta'$. We set $S'_x=\bigcup _{i=1}^m
S'_x\Delta_i$, $S'=\bigcup _{x\in K'} S'_x$ and $S'\Delta=\bigcup
_{x\in \partial \Delta\cap K'} S'_x\Delta$.

For $v\in S'_x\Delta$ denote by $\theta (v)$ the angle between $v$
 and the interior normal $\nu_{\Delta}(x)$ of $\Delta'$ with respect
to $\Delta$ at $x$. Let $dx$ be the volume element on $K'$ and let
$\lambda_x$ be the Lebesgue measure on $S'_x$. We define the {\it
Liouville} measure on $S'$ by  $d\mu' (x,v)=\cos\theta(v)
d\lambda_x(v)\otimes dx$. Note that $d\mu' (x,v)\otimes dt$ is the
ordinary Liouville measure invariant under the geodesic flow on
each $n$-simplex $\Delta$ of $K$.

Therefore, for $\mu'$-$a.e.$ $v\in S'\Delta $, the geodesic
$\gamma_v$ in $\Delta$ determined by $\dot \gamma_v(0)=v$ meets
${\partial\Delta\cap K^{(n-1)}\setminus K^{(n-2)}}$ after a finite
time $t_v>0$ so that $I(v)=-\dot\gamma_v(t_v) \in S'\Delta$. Note
that $\gamma_v(t_v) \in K'$ since $K$ is boundaryless. $\mu'$ is
invariant under the involution $I$.

Let $I(v)=u+\cos\theta(I(v)) \nu _{\Delta_n}(\gamma_v(t_v))$,
where $u$ is tangent to $K'$ and set $F(v)=\bigcup_i
\{-u+\cos\theta(I(v)) \nu_ {\Delta_n^i}(\gamma_v(t_v))\}$, where
the union is taken over all $n$-simplices containing
$\gamma_v(t_v)$ except $\Delta$.

Thus there is a subset $S_1\subset S'$ of full $\mu'$-measure such
that $F(v)$ is defined for any $v\in S_1$. We set recursively
$S_{i+1}=\{ (x,v)\in S_1 \backslash F(v)\subset S_i  \}$ and
define  $S_\infty =\bigcap_{i=0}^\infty S_i$, $V=S_\infty \cap
I(S_\infty)$. By construction, $V$ has full $\mu'$-measure.

We define the geodesic flow on the space $SK$ (or $TK$) in the
following way:
$$\text{For $(x,v)\in V$, put }\cases g^t(x,v) & = (X_{(x,v)}(t),\dot X _{(x,v)}(t))
       \\  g^0(x,v)&=(x,v),\endcases $$ where $g^t$ is the
       ordinary geodesic flow in the interior of every $n$-simplex
       and in the case where , for $t_0\in \Bbb R^+$
       $X_{(x,v)}(t_0)\in K'$, we set $\dot X_{(x,v)}(t_0)= \dot
       X_{(x,v)}(t_0+)$ (therefore, $\dot X_{(x,v)}(t_0)\in F(\dot X_{(x,v)}(t_0-)$).

 \subhead { 1.4 General Markov process}\endsubhead

Assume that $K$ is an admissible $n$-dimensional Riemannian
complex, with the metric $g$ and corresponding distance function
$d$. When $K$ is not compact, let $K_D= K\bigcup \{D\}$ be the
one-point compactification of $K$. Then, we can define a metric
$\delta$ on $K_D$ such that the topology on $K$ generated by
$\delta$ is the same as the topology generated by $d$. In case $K$
is already compact, we simply adjoint $D$ as an isolated point and
define the metric $\delta$ on $K_D$ by letting $d=\delta$ on
$K\times K$ and $\delta(p,D)=1$ for $p\in K$. Therefore, the
restriction of $\delta$ to $K\times K$ is uniformly continuous
with respect to $d$.

Denote by $C(K)$ the space of bounded continuous real-valued
functions on $K$, $C_0(K)$ the subspace of $C(K)$ such that the
functions have a null limit at infinity  and $C_c(K)$ the space of
functions in $C(K)$ with compact support. Clearly, these tree
spaces are the same if $K$ is compact. $C(K)$ endowed with supnorm
is a (real) Banach space and $C_0(K)$, $C_c(K)$ are Banach
subspaces of $C(K)$. The space $C_c(K)$ is dense in the space
$C_0(K)$.

Finally, whenever the term measurable is used it will refer to the
basic $\sigma$-algebra of Borrel sets in $K$ (or $K_D$).

 The usual setup for the theory of {\it temporarily homogeneous
 Markov process} defined on measurable space $(\Omega\times [0,\infty[, \frak M\times \frak
 R)$ ($\frak R$ is the Borrel $\sigma$-algebra in
 $[0,\infty[$) with values in topological measurable space $(E,\frak B)$ is to
 consider the following objects :
 \roster
 \item We adjoint a point $D$ to the space $E$. We write
 $E_D = E\cup \{D\}$ and $\frak B_D$ the $\sigma$-algebra in $E_D$
 generated by $\frak B$.

 \item For each $x\in E_D$, a probability measure $P_x$ on
 $(\Omega,\frak M)$.

\item An increasing family (a {\it filtration})$(\frak M_t)_{t\ge
0}$ of sub-$\sigma$-algebras of $\frak M$ and distinguished point
$\omega_D$ of $\Omega$.

\item For each $t\in [0,\infty[$ a measurable map $Y_t : (\Omega,
\frak M) \rightarrow ( E_D,\frak B_D)$ such that if
$Y_t(\omega)=D$ then $Y_s(\omega)=D$ for all $s\ge t$,
$Y_\infty(\omega) =D$ for all $\omega$ and $Y_0(\omega_D)= D$.

\item For each $t\in [0,\infty[$ a {\it translation operator}
$\theta_t: \Omega \rightarrow \Omega$ such that
$\theta_\infty\omega = \omega_D$ for all $\omega$.

\endroster

We call the collection $Y=(\Omega,\frak M,\frak M_t,
Y_t,\theta_t,P_x)$ a (temporally homogeneous) {\it Markov process}
with state space $(E,\frak B)$ if and only if the following axioms
hold :

 \roster

  \item For each $t\ge 0$ and fixed $\Gamma \in \frak B$, the
  function $x\mapsto P(t,x,\Gamma)=P_x\{Y_t \in \Gamma\}$ is $\frak B$ measurable.
\item For all $x\in E$, $ P(0,x,E\setminus\{x\})=0$ and
$P_D\{X_0=D\}=1$.

\item For all $t,h\ge0$, $Y_t\circ\theta_h=Y_{t+h}$ (homogeneity).

\item For all $s$, $t\in \Bbb R^+$, $x\in E_D$ and $\Gamma \in
\frak B_D$,
 $P_x\{X_{t+s}\in \Gamma|\frak M_t\}= P(s,X_t,\Gamma)$ (Markov property).

 \endroster

The point $D$ may be always thought of as a "cemetery" when we
regard $t\mapsto Y_t(\omega)$ as the trajectory of particle moving
randomly in the space $E$. With this interpretation in mind, we
name the random variable  $\xi(\omega)=\inf\{t; X_t(\omega)=D\}$
the {\it lifetime}.

\head { 2. Isotropic transport process}\endhead

In this section, $K$ will denote a complete  admissible Riemannian
complex with dimension $n$ and we will use all notations of the
first section.

\subhead{ 2.1 An intuitive approach}\endsubhead

Let $\Sigma K$ denote the space of links of the complex $K$.
Choose a point $(x_0,v_0)$ from the space $\Sigma K$ and assume
that the point $x_0$ is in the topological interior of a maximal
simplex $\Delta_0$. Intuitively, a particle starting from the
point $x_0$ travels geodesically, in direction $v_0$ chosen
randomly, during exponentially distributed waiting time $s_1$ to a
new position $x_1$ supposed in the interior of $\Delta_0$. At
$x_1$, the particle chooses a new direction $v_1$ in the link
$S_{x_{1}}$ over $x_1$ with the uniform probability $\Bbb P[v_{1}
\in d\lambda]=\lambda_{x_{1}}(d\lambda)$, where $\lambda$ denotes
the normalized Lebesgue measure on $S_{x_{1}}$. From the point
$x_1$ and in the direction $v_1$, the particle travels
geodesically during exponentially distributed waiting time $s_2$
to a position $x_2$ in the interior of the simplex $\Delta_0$. So
the particle continues its motion in the interior of $\Delta_0$
until it hits transversally (because of the construction of the
generalized geodesic flow on the admissible complexes) the border
of the simplex $\Delta_0$ at an interior point of a
$(n-1)$-simplex adjacent to $\Delta_0$. Note this hit point $x_n$.
Starting now from $x_n$ and choosing randomly a new direction in
the link over $x_n$, the particle travels geodesically during
exponentially waiting time $s_n$ to a new position in the interior
of a maximal simplex (which could be $\Delta_0$) and so on.

\subhead{2.2 Mathematical approach}\endsubhead

Right now, we will give a mathematical form to the random walk
just described above.

Consider the product space $L = \Sigma K\times \Bbb R^+$ and the
product $\sigma$-algebra $\frak F =
 \frak E \times \frak B$, where $\frak E$ and $\frak B$ are
 respectively Borrel $\sigma$-algebra  of $\Sigma K$ and of
 $\Bbb R^+$. Note $\Omega = L^{\Bbb N}$ and $ \frak G ={\frak F}^{\Bbb
 N}$, where $\Bbb N$ is the set of positive entire numbers. Thus $(L,\frak F)$
 and $(\Omega, \frak G)$  are measurable spaces and the points $\omega \in \Omega$ are
 sequences  $\{((x_l,v_l),t_l)\in \Sigma K\times \Bbb
 R^+ ; l\in \Bbb N\}$.

 Let $\{((x_l,v_l),t_l)\in \Sigma K\times \Bbb
 R^+ ; l\in \Bbb N\}$ be a point of $\Omega$ and set
 $\tilde Y_{k}(\omega)=((x_{k},v_{k}),t_{k})$,
 $Z_{k}(\omega)=(x_{k},v_{k})$ and $\tau_{k}(\omega)=t_{k}$.
The functions $\tilde Y_{k} : (\Omega, \frak G)\flech (L,\frak
 F)$, $Z_{k} : (\Omega, \frak G)\flech (\Sigma K, \frak E)$ and
  $\tau_{k} :(\Omega, \frak G)\flech (\Bbb R^+,\frak B)$ are
  measurable.

  Finally, we shall consider the following space of events:
  $$ \Omega'=\{\omega \in \Omega \mid \forall k\in \Bbb N, Z_{k+1}(\omega)\neq
  Z_{k}(\omega), \tau_0=0, \tau_{k+1}(\omega)>
  \tau_{k}(\omega)\}.$$

  Put $\xi(\omega) = \lim_{n\rightarrow\infty}\tau_{n}(\omega)$
  (life time) and let $K_D= K\bigcup \{D\}$ denote
  the one point compactification. The space $K$ is assumed
  semi-compact so we shall endow $K_D$ with a metric
  $d'$ such that the space $(K_D,d')$ is compact and the
  restriction of $d'$ to $K$ coincides with the beginning  metric of $K$.

  Now, we will define the $K$-valued geodesically random walk which interests us.
  Let, for $t\geq 0$ :

$$Y_t(\omega) =   \cases  X_{Z_i(\omega)}(t-\tau_i(\omega))&\text{if $\tau_i(\omega) \le t\le
\tau_{i+1}(\omega)$ ,}\\ D&\text{if $\xi(\omega)\leq t$
 ,}\endcases$$
 where $X$ is the $K$-projection of the generalized geodesic flow on the complex $K$.
  According to the latest definition, we have for every $\omega \in \Omega$, $Y_{\infty}=D$.

\subhead{ 2.3 Markov property}\endsubhead

In the following paragraph, we will  complete the preceding
construction to define the admissible complex-valued {\it
isotropic transport process} and then we will show that the last
process is a strong Markov one.

Let $K$ denote an admissible Riemannian complex and define the
next transition density on the measurable space $(L,\frak F)$ as :

$$N(z,t;dz,ds) =   \cases 0 &\text{if $t<s$ ,}\\
  \lambda_{x}(dz) e^{-(s-t)} ds &\text{ if $s\geq t $
 ,}\endcases$$

with $z=(x,v)$, $dz=(x,dv)$ and $\lambda_{x}$ is the uniform
measure on the link $S_xK$.

\proclaim{ Proposition 2.1}

Let $\gamma$ denote a probability measure on the  measurable space
$(L,\frak F)$. Then, there exists a probability measure $P^\gamma$
on the measurable space $(\Omega,\frak G)$ such that the
 coordinate mappings $\{ \tilde Y_n ; n\in \Bbb N\}$ form a
 temporally homogeneous Markov process on the measured space
 $(\Omega,\frak G,P^\gamma)$, with $\gamma$ as initial distribution
and $N$ the transition function, i.e:
 $$P^\gamma (\tilde Y_{n+1}\in A | \tilde Y_{0},\ldots,\tilde
 Y_{n}) = \int_{A }N(Z_{n},\tau_{n};dz,ds),$$
for all $A$ belonging to $\frak F$ and $n\in \Bbb N$.
\endproclaim

\demo{Proof}

The proposition is an immediate corollary  of {\it I.Tulcea's}
theorem (see [12] pp. 613-615). \hfill $\square$
 \enddemo

 If $\gamma $ is the measure $\lambda_{x}\otimes \delta_{0}$,
 with $\delta_{0}$ the Dirac mass at $0\in \Bbb R$, then we
 will write $P^{\lambda_{x}}$ or $P^{x}$ for $P^\gamma $.
 Consequently, we have for every $x\in K$,
 $P^x (\Omega')=1$ and the process $\{ \tilde Y_n ; n\in \Bbb N\}$ will be
 Markov on the measured space $(\Omega',\frak
 G',P^x)$. We will note $\Omega$ the set of sequences
 $\{(z_n,t_n)\in L; n\ge 0\}$ where
 $z_{n+1}\neq z_n$ and $0=t_0<t_1<\ldots <t_n<\ldots$,
 and $\frak G$ will denote the $\sigma$-algebra of $\Omega$
 generated by $\{ \tilde Y_n ; n\in \Bbb N \}$. Thus, we will use
 in the following, the probability space(s) $(\Omega,\frak G,P^x)$.

Right now, let $(Y_t)_{t\ge 0}$ denote the $K$-valued random walk
constructed in the last section (3.1). For all $\omega\in \Omega$,
the map $t\mapsto Y_t(\omega)$ is continuous on $\Bbb R^+$ and has
 left-hand limits on $[0,\xi (\omega )[$. We complete the
$\sigma$-algebra $\frak G$ by adjoining a point $\omega_D$ to
 $\Omega$ with $Y_t(\omega_D)=D$ for all $t$,
 $\{\omega_D\} \in \frak G$ and $P^x(\{\omega_D \})=0$
for all $x\in K$. We set $Z_n(\omega_D)= D$ and
 $\tau_n(\omega_D)=\infty$ for all $n\in \Bbb N$ and note $P^D$
 the Dirac mass at $\omega_D$.

 Next we define the translation operators  $(\theta_t)_{t\ge0}$ as
 follows:
   for all $t\ge 0$, $\theta_t\omega_D=\omega_D$ ; if
   $t\ge \xi(\omega)$ then $\theta_t\omega=\omega_D$, while if
   $t_k\le t <t_{k+1}$, $k\ge 0$ then $\theta_t\omega =
   \{(z_{n+k},(t_{n+k}-t)\vee 0); n\ge 0\}$, where $\omega= \{ (z_n,t_n) ; n\ge
   0\}$.

   Thus, we have $Y_s\circ \theta_t
   =Y_{s+t}$ for all $s,t \in \Bbb R^+$.

   \definition{Definition 2.2}

   We call the stochastic process $Y= (\Omega,\frak
   G,Y_t,\theta_t,P^x)$ the (an) isotropic transport process (motion)
    with values in the admissible Riemannian complex $K$.
   \enddefinition

  Let $\frak G_n:= \sigma \{\tilde Y_i ; 0\le i\le n\}$
   and $\frak F^0_t:= \sigma \{ Y_s ; s\le t\}$ denote
   respectively the $\sigma$-algebra of $\Omega$ generated by $\{\tilde Y_i ; 0\le i\le
    n\}$ and the one generated by $\{ Y_s ; s\le t\}$.

    \proclaim{Lemma 2.3}

    Let $\Lambda \in \frak F^0_t$ ; then, for all $n\ge0$, there
    exists $\Lambda_n \in \frak G_n$ such that :
    $$\Lambda\cap \{\tau_n \le t <\tau_{n+1}\}=
    \Lambda_n\cap\{t< \tau_{n+1}\}.$$

    \endproclaim

  \demo{Proof}

  Note : $$\frak G_t := \sigma \{\Lambda \in \frak F^0_t |(\forall
  n\ge 0) (\exists \Lambda_n \in \frak G_n),\Lambda\cap \{\tau_n \le t <\tau_{n+1}\}=
    \Lambda_n\cap\{t< \tau_{n+1}\}\}.$$
    We can easily check that, for all $A\in \frak
    E_D$, the sets $\{Y_s\in A\}_{s\le t}$ belong
    to the $\sigma$-algebra $\frak G_t$. Thus, we end the proof; indeed the sets $\{Y_s\in A\}_{s\le t}$
    generate the $\sigma$-algebra $\frak F^0_t$.

\hfill $\square$
  \enddemo

  We set, for real functions  $g\in C_0(\Sigma K)$ and $f\in
  C_0(K)$ (or simply measurable functions ):
   \roster

  \item $Pg(x):=\int _{\Sigma_xK}g(x,\eta) d\lambda_x(\eta)$ .

  \item $\forall t>0 , T_t^0f(x) := \int _{\Sigma_xK}f(X_{(x,\eta)}(t))
  d\lambda_x(\eta)$ and $T_tf(x) := E^x[f(Y_t)]$, the latest is the expectation with respect to $Y_t$.

   \item $\forall \lambda >0, R^0_\lambda f(x):= \int_{\Bbb R^+} e^{-\lambda
  t}T_t^0f(x)dt$ and $R_\lambda f(x) := \int_{\Bbb R^+} e^{-\lambda
  t}$ $T_t f(x)dt$, respectively the resolvent operator of $T_t^0$ and of $T_t$.

\endroster

 \proclaim{Proposition 2.4}

Let $f\in C^0(K)$, then, for all $\lambda >0$ we have:

$$R_\lambda f = \sum^ \infty _{n=0}(R^0_{\lambda+1})^{n+1}f,$$
where $(R^0_{n+1})^0 := Id $ the identity map.
  \endproclaim

\demo{Proof}

 First we write:

$$R_\lambda f(x) = [\int_0^{\tau_1} + \sum_{i=1}^{\infty}\int_{\tau_i}^{\tau_{i+1}}]e^{-\lambda
  t} f(Y_t)dt .$$

  Taking into account the distribution  of $\tau_1$ and the initial distribution
 of the process $Y$, the first integral becomes:

  $$\int_0^\infty e^{-(1+\lambda)s}T_s^0f(x) ds = R^0_{1+\lambda}f(x).$$

 For the second part of the decomposition, we will prove by
induction argument that for all $i\ge 1$ the following equality:

$$\circledR \hskip 1cm [\int_{\tau_i}^{\tau_{i+1}}e^{-\lambda  t} f(Y_t)dt] =
(R^0_{\lambda+1})^{i+1} f(x).$$

Let see the case $i=1$:

$$[\int_{\tau_1}^{\tau_{2}}e^{-\lambda  t} f(Y_t)dt] = [e^{-\lambda
\tau_1 }\int_{0}^{\tau_{2}-\tau_1}e^{-\lambda  t}
f(Y_{t+\tau_1})dt] , $$
 which is equal to :

$$[e^{-\lambda \tau_1 }(R^0_{\lambda+1})f(X_{Z_1}(0) )] =
[e^{-\lambda \tau_1 }(PR^0_{\lambda+1})f(X_{Z_0}(\tau_1) )].$$

Using the distribution of $\tau_1$, we obtain:

$$ (R^0_{\lambda+1}) (R^0_{\lambda+1})f(x).$$

Assume the property $\circledR$ until the order $l$, and see what
happens at the order $l+1$:

$$[\int_{\tau_{l+1}}^{\tau_{l+2}}e^{-\lambda  t} f(Y_t)dt] = [e^{-\lambda
\tau_{l+1} }\int_{0}^{\tau_{l+2}-\tau_{l+1}}e^{-\lambda  t}
f(Y_{t+\tau_{l+1}})dt] , $$ which is equal to:
$$[e^{-\lambda \tau_{l+1} }(R^0_{\lambda+1})f(X_{Z_{l+1}}(0) )] =
[e^{-\lambda \tau_l}
e^{-\lambda(\tau_{l+1}-\tau_l)}(PR^0_{\lambda+1})f(X_{Z_l}(\tau_{l+1}-\tau_l)
)].$$

Using the distribution of  $(\tau_{l+1}-\tau_l) $:

$$ =[e^{-\lambda \tau_l} (R^0_{\lambda+1}) (R^0_{\lambda+1})f(X_{Z_l}(0))],$$
this latest expectation is equal to:

$$[\int_{\tau_{l}}^{\tau_{l+1}}e^{-\lambda  t}
R^0_{\lambda+1}f(Y_t)dt].$$

 Hence, using the recurrence hypothesis applied to the
function $R^0_{\lambda+1}f$, we obtain the equality $\circledR$ at
the order $l+1$.

For the end of the proof, note that the series $\sum^ \infty
_{n=0}(R^0_{\lambda+1})^{n+1}f$ converges uniformly because we
have for all function $f\in C^0(K)$, the estimation
$||R^0_{\lambda+1}|| \le \frac{1}{\lambda+1}$ (the $\sup$ norm),
which end the proof. \hfill $\square$
 \enddemo

\proclaim{Lemma 2.5}

 Let $f$ be a measurable real (positive) function on $(K,\frak
 B)$. Then, for all $t\ge 0$ and $\lambda >0$, we have:

 $$E\{\int_t^\infty e^{-\lambda u} f(Y_u)du | \frak F_t^0 \} = e^{-\lambda t} R_\lambda f(Y_t). $$

\endproclaim

\proclaim{Remark 2.6}

By the lemma 2.3 of this paragraph, to establish the lemma 2.5 it
suffices to show the  same equality(s) on the sets
 $\Lambda_n \in \frak F^0_t$ with
    $$\Lambda_n\cap \{\tau_n \le t <\tau_{n+1}\}=
    \Lambda_n\cap\{t< \tau_{n+1}\}.$$

i.e:

$$\maltese \hskip 1cm E\{\int_t^\infty e^{-\lambda u} f(Y_u)du | \Lambda_n \}
= E\{ e^{-\lambda t} R_\lambda f(Y_t) | \Lambda_n \}. $$

\endproclaim

\demo{Proof of Lemma 2.5}

Consider the left side of the equality $\maltese$ and set it in
the following way :
$$E\{\int_t^\infty e^{-\lambda u} f(Y_u)du | \Lambda_n \}
= [(\int_t^{\tau_{n+1}} +
\sum_{i=n+1}^{\infty}\int_{\tau_i}^{\tau_{i+1}})e^{-\lambda
  u} f(Y_u)du | \Lambda_n].$$

  Using, the Markov property  of the process $\{\tilde Y_n; n\ge
  0\}$, the fact that $\Lambda_n \subset \{ \tau_n \le t\le
  \tau_{n+1}\}$ and the  exponential distribution of  the random variable $\tau_{n+1}-t\wedge
  \tau_{n+1}-\tau_n$, the first integral of the decomposition becomes:
  $$[ e^{-\lambda t}e^{-( t-\tau_n)}\int_0^\infty e^{-(\lambda+1)u}Pf(X_{Z_n}(u+(t-\tau_n))) du |
  \Lambda_n]$$
  which is equal to:
  $$[e^{-\lambda t} e^{-( t-\tau_n)}R^0_{\lambda+1}f(X_{Z_n}(t-\tau_n)) | \Lambda_n].$$

For the second half of the decomposition, we will show by
induction argument that for all $i\ge 1$, we have the equality:
$$\circledR_n^t \hskip 1cm [\int_{\tau_{n+i}}^{\tau_{n+i+1}}e^{-\lambda  u} f(Y_u)du | \Lambda_n] =
[e^{-\lambda t}e^{-( t-\tau_n)}(R^0_{\lambda+1})^{i+1}
f(X_{Z_n}(t-\tau_n)) | \Lambda_n]. $$

Let see the case $i=1$:
$$[\int_{\tau_{n+1}}^{\tau_{n+2}}e^{-\lambda  u} f(Y_u)du | \Lambda_n] = [e^{-\lambda
\tau_{n+1} }\int_{0}^{\tau_{n+2}-\tau_{n+1}}e^{-\lambda  u}
f(Y_{u+\tau_{n+1}})du | \Lambda_n] , $$ which is equal to:

$$[e^{-\lambda \tau_{n+1}}(R^0_{\lambda+1})f(X_{Z_{n+1}}(0) ) | \Lambda_n]=
[e^{-\lambda t} e^{-\lambda
(\tau_{n+1}-t)}(R^0_{\lambda+1})f(X_{Z_{n}}(\tau_{n+1}-\tau_{n}) )
| \Lambda_n],$$ which is the same as:
$$[e^{-\lambda t} e^{-\lambda
(\tau_{n+1}-t)}(R^0_{\lambda+1})f(X_{Z_{n}}((\tau_{n+1}-t)+(t-\tau_{n}))
) | \Lambda_n].$$

 Using the Markov property of $\{\tilde Y_n; n\ge  0\}$ and the distribution
 of $(\tau_{n+1}-t)\wedge(\tau_{n+1}-\tau_{n})$, we obtain :

$$[e^{-\lambda t}e^{-(t-\tau_{n})}
\int_0^\infty
e^{-(\lambda+1)u}P(R^0_{\lambda+1})f(X_{Z_n}(u+(t-\tau_{n}))) du |
\Lambda_n],$$ what is equal to:

$$ [e^{-\lambda t}e^{-(t-\tau_{n})}
R^0_{\lambda+1}(R^0_{\lambda+1})f(X_{Z_n}(t-\tau_{n})) |
\Lambda_n].$$

Now assume that the property $\circledR_n^t$ comes true until the
order $l$, and let's see what will happen at order $l+1$.

$$[\int_{\tau_{n+(l+1)}}^{\tau_{n+(l+2)}}e^{-\lambda  u} f(Y_u)du | \Lambda_n]
= [e^{-\lambda \tau_{n+(l+1)}
}\int_{0}^{\tau_{n+(l+2)}-\tau_{n+(l+1)}}e^{-\lambda  u}
f(Y_{u+\tau_{n+(l+1)}})du | \Lambda_n] , $$ what is equal to :
$$[e^{-\lambda \tau_{n+(l+1)} }(R^0_{\lambda+1})f(X_{Z_{n+(l+1)}}(0) ) | \Lambda_n] ,$$
which is equal to:
$$ [e^{-\lambda \tau_{n+l}}
e^{-\lambda(\tau_{n+(l+1)}-\tau_{n+l})}(PR^0_{\lambda+1})f(X_{Z_{n+l}}(\tau_{n+(l+1)}-\tau_{n+l})
)| \Lambda_n].$$

Then, using the distribution of  $(\tau_{n+(l+1)}-\tau_{n+l}) $ we
ge:

$$ =[e^{-\lambda \tau_{n+l}} (R^0_{\lambda+1}) (R^0_{\lambda+1})f(X_{Z_{n+l}}(0)) | \Lambda_n],$$
this latest expectation is equal to:

$$[\int_{\tau_{n+l}}^{\tau_{n+(l+1)}}e^{-\lambda  u}
R^0_{\lambda+1}f(Y_u)du | \Lambda_n].$$

 Thus, if we apply the recurrence hypothesis to the
function $R^0_{\lambda+1}f$, we obtain the equality
$\circledR_n^t$ at the order $l+1$.

 Up to now  we have shown that the left side of the equality
$\maltese$ is equal to:
$$[ e^{-\lambda t}e^{-( t-\tau_n)}
\{R^0_{\lambda+1}f(X_{Z_n}(t-\tau_n))
+\sum_{i=1}^{\infty}(R^0_{\lambda+1})^{i+1} f(X_{Z_n}(t-\tau_n))\}
| \Lambda_n].$$ Thus by  proposition 2.4 of this section, this sum
is equal to:
$$[ e^{-\lambda t}e^{-( t-\tau_n)}
R_{\lambda}f(X_{Z_n}(t-\tau_n)) | \Lambda_n].$$

 Using once again the Markov property of
$\{\tilde Y_n; n\ge  0\}$,  we get:
$$[ e^{-\lambda t}e^{-( t-\tau_n)}
R_{\lambda}f(X_{Z_n}(t-\tau_n)) | \Lambda_n]= [ e^{-\lambda t}
R_{\lambda}f(Y_t) | \Lambda_n]$$

which was to be proved. \hfill $\square$
\enddemo

Right now, we have collected all the ingredients to prove the
following theorem:

\proclaim{Theorem 2.7}

 Let $Y= (\Omega,\frak F^0_t,Y_t,\theta_t,P^x)$ be the isotropic
transport process with  values in the admissible Riemannian
complex $K$. Then $Y$ is a strong Markov process.

\endproclaim

\proclaim{Remark 2.8 ([14] I pp 97-100)}

It suffices to show that the process $Y$ is a Markov process
because the right continuous (with right continuous trajectories)
Markov process is always strongly Markov for the filtration
$\{\frak F_{t+}^0\}$. But we know that, in case of continuous
stochastic process, the filtration $\{\frak F_{t+}^0\}$ is equal
to the filtration $\{\frak F_{t}^0\}$, which includes the case of
the isotropic transport process (it is trajectories continuous).
\endproclaim

\demo{Proof of Theorem 2.7}

By Lemma 2.5 of this section we have:

$$E\{\int_t^\infty e^{-\lambda u} f(Y_u)du | \frak F_t^0 \} =
E^{Y_t} \{\int_0^\infty e^{-\lambda(t+u)}f(Y_u)du\}. $$

Then, if the function $f$ is bounded we have the next equality:

$$E\{\int_t^\infty \varphi(u) f(Y_u)du | \frak F_t^0 \} =
E^{Y_t} \{\int_0^\infty \varphi(t+u)f(Y_u)du\}, $$ whenever the
function $\varphi$  is a linear combination of exponentials and
hence, by uniform approximation, whenever $\varphi$ is continuous
and vanishes at infinity. Then, consider the following sequence of
functions:

$$\varphi_n(s+t+u) =   \cases  0 &\text{if $\frac{1}{n} \le u$ ,}
\\ \frac{1}{n}-x &\text{if $0\leq u < \frac{1}{n}$ .}
 \endcases$$

 The sequence $(\varphi_n)_{n\ge 0}$ is a sequence of continuous functions vanishing at infinity
 and converging to the Dirac mass at $s+t$, while the
 map $u\mapsto f(Y_u)$  is a bounded (right) continuous function.
 Consequently, if we take the limit we obtain:

$$E\{ f(Y_{t+s})du | \frak F_t^0 \} =
E^{Y_t} \{f(Y_s)\}. $$

In other words, $Y$ is Markov process. \hfill $\square$
\enddemo

\head { 3.Wiener measure.}
\endhead

\subhead{ 3.1 Construction}\endsubhead

Let $Y= (\Omega,\frak F^0_t,Y_t,\theta_t,P^x)$ be the isotropic
transport process in the complete admissible Riemannian complex
$K$ constructed in the last section. Set for a real $\eta>0$ and
$z=(x,v)\in \Sigma K$, $\eta z := (x,\eta v)$.

Define now a process $Y^{\eta}$ from  $Y= (\Omega,\frak
F^0_t,Y_t,\theta_t,P^x)$ in the following way :

$$Y_t^{\eta}(\omega) =   \cases  X_{\eta Z_i(\omega)}(\frac{t}{\eta^{2}}-\tau_i(\omega))
&\text{if $\tau_i(\omega) \le \frac{t}{\eta^{2}}\le\tau_{i+1}(\omega)$ ,}\\
D &\text{ if $\xi(\omega)\leq \frac{t}{\eta^{2}}$
 .}\endcases$$

Thus the process $Y^\eta= (\Omega,\frak
 F^0_t,Y^\eta_t,\theta_t,P^x)$ is (trajectories) continuous and it is, as the process $Y$,
 strongly Markov.

 \proclaim{ Proposition 3.1}

Let $K$ be an admissible Riemannian and $C(\Bbb R^+, K)$ be the
space of continuous paths in $K$. Then for each $\eta>0$, the
process $Y^\eta$ generates a measure $\mu_{\eta}$ on the space
$C(\Bbb R^+, K)$.
 \endproclaim

\demo{Proof}

For $\eta >0$, set $P^{\eta}_{s,t}(p,A)$ with $p\in K$ and $A\in
\frak B(K_{D})$, the transition probability of the process
$Y^{\eta}$ (i.e. $P^{\eta}_{s,t}(p,A) := Prob \{Y^{\eta}_{t+s}\in
A; Y^{\eta}_{s}=p\}$ ).

Consider the finite sets of reals $J=\{t_1<t_2<\ldots<t_n\}\subset
(\Bbb R^+)^n$. Then, for each finite set
$J=\{t_1<t_2<\ldots<t_n\}$, we define probability measure in the
following way:
$$\text{ for $B\subset K_D^n$, $ P^{\eta}_J(B)=\int_B P^x(dx_0)\int
P^{\eta}_{0,t_1}(x_0,dx_1)\int \ldots \int
P^{\eta}_{t_{n-1},t_n}(x_{n-1},dx_n)$ .}$$

Let $\Phi ( \Bbb R^+)$ denote the set of the finite subset of
$\Bbb R^+$. Then, the system $\{P^{\eta}_J; J\in \Phi ( \Bbb
R^+)\}$, and thanks to the Markov property of $Y^\eta$, is a
 projective system on $(K_D,\frak B(K_D))$ (i.e : if $\pi_J^I$
 (respectively $\pi_J$) is the natural projection of $K^I$ (respectively $\Omega$))
 to $K^J$ then $P^{\eta} _I(\pi _J ^I)^{-1}=P^{\eta} _J$).

 On the other hand,
 the trajectories of $Y^\eta$ are continuous and the space $K$ is
Hausdorff and $\sigma$-compact. Consequently, and using the {\it
Kolmogorov} theorem [5], we get a probability measure $\mu_{\eta}$
on the space $C(\Bbb R^+,K)$.\hfill $\square$
\enddemo

\subhead{ 3.2. Wiener measure} \endsubhead

Right now, we will announce the main theorem of this section.

\proclaim{ Theorem 3.2}

Let $K$ denote an admissible Riemannian complex, and consider the
family  $\{Y^\eta \}_{\eta
>0 }$ of the isotropic transport processes constructed in the paragraph above
and let $(\mu_{\eta})_{\eta >0}$ be the family of the generated
probability measures on $C(\Bbb R^+,K)$. The space $C(\Bbb R^+,K)$
is provided with the compact-open topology. Then the family
$(\mu_{\eta})_{\eta >0}$ has a convergent subsequence.
\endproclaim

To prove the last theorem we need the following lemma:

\proclaim{Lemma 3.3}

Under the hypothesis of Theorem 3.2, the family of the probability
measures $(\mu_{\eta})_{\eta
>0}$ is Tight, i.e :

$$\lim \Sb {\eta \rightarrow 0}\\{c\rightarrow 0}\endSb
Prob\{ \sup \Sb {t-c<t_1<t_2<t+c}\\{0\le t_1<t_2\le N} \endSb
\min[d(Y^\eta_{t_1},Y^\eta_{t});d(Y^\eta_{t},Y^\eta_{t_2})]>\epsilon
\} = 0 .$$
\endproclaim

\proclaim{Remark 3.4}

Before proceeding to look at the proof of the lemma, recall fist
the two following facts:

\roster \item When the space $C(\Bbb R^+,K)$ is provided with the
compact-open topology, the Tightness property is equivalent,
following Stone [29], to the equality of the Lemma 3.3.

\item According to J\o rgensen's result [19, Lemma 1.4], if the
following property :
$$\text{for all $\epsilon>0$, there exists $\alpha>0$ such that,
   $\sup \Sb {p\in K_D}\\ 0<s \endSb \frac{1}{s}
   P^{\eta}_{0,s}(p,B_D^c(p,\epsilon))\le \alpha $ ,}$$
   comes true then the equality of Lemma 3.3 is also true.

\endroster

\endproclaim

\demo{Proof of Lemma 3.3}

By Remark 3.4, if we show the following :
$$\text{$\forall \epsilon>0$, $\exists \alpha>0$ ,
   $\lim_{ \eta\rightarrow 0} \sup \Sb {p\in K_D}\\ 0<t \endSb
\frac {Prob\{Y_t^{\eta}\in B_D^c(p,\epsilon)\}}{
\frac{t}{\eta^2}}\le \alpha $ ,}$$ then the sequence
$(\mu_{\eta})_{\eta >0}$ is Tight.

We will assume that $\epsilon <\eta $ (otherwise, the probability
needed should be null) and $\frac{t}{\eta^{2}}< \tau_1$ (see the
recurrence in proof of Lemma 2.5) which doesn't affect the result.
On the other hand, $\epsilon$ is necessarily lower or equal than
$\frac{t}{\eta^2}$ unless the sought after probability should
vanish and then, there is nothing to prove.

 Thus, we have:
$$ Prob\{Y_t^{\eta}\in B_D^c(p,\epsilon)\}=
 E\{ I_{B_D^c(p,\epsilon)}(Y_t^{\eta}) | \epsilon \le \frac{t}{\eta^{2}}<
\tau_1 \}.$$ Using the Markov property, we obtain:
 $$Prob\{Y_t^{\eta}\in B_D^c(p,\epsilon)\}=E\{ e^{-\frac{t}{\eta^2}}\int_{0}^{\infty}P I_{B_D^c(p,\epsilon)}
 (X_{(p,\eta \zeta)}(\frac{t}{\eta^{2}}+s)) e^{-s} ds | \epsilon \le \frac{t}{\eta^{2}} \},$$
that is equal to:
 $$ E\{ e^{-\frac{t}{\eta^2}} R^0_1 I_{B_D^c(p,\epsilon)}
 (X_{(p,\eta \zeta)}(\frac{t}{\eta^{2}}))|\epsilon \le \frac{t}{\eta^{2}} \} .$$

 Using the fact that $||R^0_1|| \le 1$ we obtain the following estimation :

 $$ E\{ e^{-\frac{t}{\eta^2}} R^0_1 I_{B_D^c(p,\epsilon)}
 (X_{(p,\eta \zeta)}(\frac{t}{\eta^{2}}))|\epsilon \le \frac{t}{\eta^{2}} \}\le
  e^{-\frac{t}{\eta^2}} .$$

So for all $t>0$ we get :

$$\frac {Prob\{Y_t^{\eta}\in B_D^c(p,\epsilon)\}}{
\frac{t}{\eta^2}} \le
\frac{e^{-\frac{t}{\eta^2}}}{\frac{t}{\eta^2}} .$$

Thus, for all $t>0$, if $\eta$ goes to zero,
$\frac{e^{-\frac{t}{\eta^2}}}{\frac{t}{\eta^2}}$ goes also to
zero, which was to be proved.

\hfill $\square$
\enddemo

Right now, we are ready to prove Theorem 3.2.

\demo{Proof of Theorem 3.2}

Consider the space $C(\Bbb R^+, K)$ provided with the compact-open
topology, where $K$ is an admissible Riemannian complex. Let
$(\mu_\eta)_{\eta>0}$ be the sequence of probability measures
generated by the family of isotropic transport processes
$\{Y^{\eta}\}_{\eta>0}$.

By Lemma 3.3, the sequence $(\mu_\eta)_{\eta>0}$  is Tight;
moreover, the space $C(\Bbb R^+, K)$ endowed with the compact-open
topology is separable. Thus, using {\it Prohorov}'s theorem (see
[5]), the sequence $(\mu_\eta)_{\eta>0}$ is relatively compact.
The proof is now complete. \hfill $\square$
\enddemo

We showed above that the sequence $(\mu_\eta)_{\eta>0}$ has a
subsequence  which converges to a probability measure. Let $W$
denote this limit; then we set the following definition :

\definition{ Definition 3.5}
The measure $W$ on the space $C(\Bbb R^+,K)$ is called a Wiener
measure.
\enddefinition

\subhead{ Example 3.6 : The smooth case}\endsubhead

Assume $K$ is a smooth Riemannian manifold of dimension $n$ and
let $\triangle $ denote the operator of Laplace-Beltrami on $K$,
then $\triangle$ is the infinitesimal generator of a Markov
process, named the Brownian motion [16], and  note it
$\{B_t^x\}_{t<\zeta'}$. Let $(U_t)_{t>0}$ denote the semigroup
associated to the Brownian motion. Suppose that, for all $f\in
C_0(K)$, $U_tf\in C_0(K)$. Then we have the following theorem:

\proclaim{ Theorem}

The sequence of processes $\{Y^{\eta}\}_{\eta>0}$ converge weakly
to the process $\{B_t^x\}_{t<\zeta'}$.
\endproclaim

\demo{Proof}

Set $T_t^{\eta}f(x)= E^x[f(Y_t^{\eta})]$ ; following a result of
{\it Pinsky} (see [27]), we have:

$$\text{ $\forall f\in C_0(K)$, $\lim_{\eta\rightarrow 0}T_t^{\eta}f =U_{\frac{t}{n}}f$,
where $n$ is the dimension of $K$.} $$

By Theorem 3.2, there exists a subsequence
$(\mu_{\eta'})_{\eta'>0}$ of the sequence of probability measures
$(\mu_\eta)_{\eta>0}$, such that $(\mu_{\eta'})_{\eta'>0}$
converges to a probability measure $W$ on the space $C(\Bbb
R^+,K)$.

Thus, by {\it Stone}'s theorem [29], $W$ is then the classical
Wiener measure generated by the Brownian motion
$\{B_t^x\}_{t<\zeta'}$. \hfill $\square$

\enddemo

\head { 4. Brownian motion.}
\endhead

By $K$ we always denote a complete admissible Riemannian complex,
consider $\{Y^\eta \}_{\eta >0 }$  the family of the isotropic
transport processes and $(\mu_{\eta})_{\eta >0}$ the corresponding
sequence of probability measures.

Let $(\mu_{\eta_k})_k$ be a subsequence of the sequence
$(\mu_{\eta})_{\eta >0}$ which converges to the Wiener measure
$W$.

Note, for $\eta_k>0$ and for each finite set
$J=\{t_1<t_2<\ldots<t_n\}$, $ P^{\eta_k}_J$ the probability
measure defined on the product space $K^n$, as follows:
$$\text{ for $B\subset K_D^n$, $ P^{\eta_k}_J(B)=\int_B P^x(dx_0)\int
P^{\eta_k}_{0,t_1}(x_0,dx_1)\int \ldots \int
P^{\eta_k}_{t_{n-1},t_n}(x_{n-1},dx_n)$.}$$

\proclaim{Proposition 4.1}

By $\Phi ( \Bbb R^+)$ we note the set of all finite subsets of
$\Bbb R^+$. Then, for all $J$ in the set $\Phi ( \Bbb R^+)$, the
sequence of probability measures $(P^{\eta_k}_J)_k$ has a
subsequence converging to a probability measure $\mu_J$ on the
space $K_D^{|J|}$ ($|J|$ is the cardinal of $J$). Moreover, the
system $\{\mu_J; J\in \Phi ( \Bbb R^+)\}$ is projective on the
space $(K_D,\frak B(K_D))$.

\endproclaim

\demo{Proof}

 Recall that for all $s\in \Bbb R^+$, $t\in \Bbb R^+$ and all $p\in K$,
 the sequence of transition functions $(P^{\eta_k}_{s,t}(p,.) )_k$ ( $P^{\eta_k}_{s,t}(p,A):= Prob
\{Y^{\eta_k}_{t+s}\in A; Y^{\eta_k}_{s}=p\}$ where  $A\in \frak
B(K_{D})$) defines a sequence of probability measures on the space
$(K_D,\frak B(K_D))$.

Moreover, the space $K_D$ is $\sigma$-compact; Thus, following
{Prohorov}'s theorem [5], there exists a probability measure
$\mu_{s,t}^{p}$ and a subsequence $(P^{\eta_k}_{s,t}(p,.) )_k$
converging weakly to $\mu_{s,t}^{p}$.

By a diagonal argument, we obtain, for all
$J=\{t_1<t_2<\ldots<t_n\}$ in $\Phi ( \Bbb R^+)$), a probability
 measure $\mu_{J}$ on  the product space $K_D^{|J|}$ in the following way:
$$\text{ for $B\subset K_D^{|J|}$, $ \mu_J(B)=\int_B P^x(dx_0)\int
\mu^{x_0}_{0,t_1}(dx_1)\int \ldots \int
\mu^{x_{n-1}}_{t_{n-1},t_n}(dx_n)$,}$$ consequently the proof is
now complete.\hfill $\square$
\enddemo

\proclaim{Remark 4.2}

The sequence $(\mu_{\eta_k})_k$ is weakly convergent to the
 Wiener measure $W$. Thus, for every set $J$ belonging to $\Phi ( \Bbb R^+)$,
 the finite dimensional distribution $W(\pi_J)^{-1}$
coincides with $\mu_J$. In particular, for all $s>0$, $t>0$ and
$p\in K_D$, we have $W^p(\pi_{\{s,t\}})^{-1}:=\mu_{s,t}^{p}$.
\endproclaim

\proclaim{Corollary 4.3}

 The function which maps a point $(t,p,\Gamma)\in \Bbb R^+\times
 K_D\times \frak B(K_D)$ to
$W(t,p,\Gamma):=W^p(\pi_{\{0,t\}})^{-1}(\Gamma)$ is a transition
function on the measurable space $(K_D,\frak B(K_D))$.
\endproclaim

\demo{Proof}

The corollary is an immediate consequence of proposition 4.1 and
remark 4.2 .\hfill $\square$
\enddemo

Just now, we are ready to give the main theorem of this section.

\proclaim{Theorem 4.4}

Let $(t,p,\Gamma)\mapsto W(t,p,\Gamma)$ denote the transition
function on the measurable space $(K_D,\frak B(K_D))$,
corresponding to the Wiener measure on the space $C(\Bbb R^+, K)$
(see corollary 4.3). Then there exists a continuous $K_D$-valued
Markov process $\{B^p_t\}_{t\ge 0}$ with $W(t,p,\Gamma)$ as
transition function.
\endproclaim

Before proceeding to look at the proof, we first give the
following definition:

\definition{Definition 4.5}

The continuous $K_D$-valued Markov process $\{B^p_t\}_{t\ge 0}$,
is called a {\it Brownian motion}.

\enddefinition

\demo{Proof of theorem 4.4}

Using a corollary of the {\it Kolmogorov}'s theorem (see [14] I
page 91 theorem 3.5), if we show that the transition functions
$(t,p,\Gamma)\mapsto W(t,p,\Gamma)$ satisfy the following two
conditions, for each compact $\Gamma\subset K_D$ : \roster \item
for all $N>0$, $ \lim_{y\rightarrow\infty}\sup_{t\le
N}W(t,y,\Gamma)=0$, \item for all $\epsilon >0$,
$\lim_{t\downarrow 0}\sup_{p\in
\Gamma}\frac{1}{t}W(t,p,B_D^c(p,\epsilon))=0$,
\endroster then the conclusion of theorem 4.4 comes true.

For the first condition, consider a  compact $\Gamma\subsetneq
K_D$ (otherwise if $\Gamma=K_D$ then the first condition is
trivially satisfied). Let $(\mu_{\eta_k})_k$ be a sequence of
measures associated to the sequence of isotropic processes which
converges weakly to the Wiener measure $W$ and let
$P^{\eta_k}_{t}(p,A):= Prob \{Y^{\eta_k}_{t}\in A;
Y^{\eta_k}_{0}=p\}$ denote the associated transition functions.

Recall that, for each $\eta>0$, all trajectories of the random
walk $Y^\eta$ are concatenations of geodesic segments, with every
geodesic segment's length lower or equal to $\eta$. Consequently,
we have for each $\eta_k>0$, $d(p,Y^{\eta_k}_t)\le \eta_k t$ if
$Y^{\eta_k}_0=p$.

Let $N>0$ some (fixed) real, then for all $t\le N$, if
$Y^{\eta_k}_0=y$, $d(y,Y^{\eta_k}_t)\le \eta_k N$. Thus, if we
consider the points $y\in K_D$ with the distance $d(y,\Gamma)$
strictly greater than $(\eta_k+\epsilon)N$, for some $\epsilon>0$,
then the probability $P^{\eta_k}_{t}(y,\Gamma)$ should vanish.

In a nutshell, we proved that, for all $\eta_k>0$, $N>0$ and $y\in
K_D$ :
$$\text{For all $\epsilon>0$, there exists $\alpha= (\eta_k+\epsilon)N$
such that if $d(y,\Gamma)\ge \alpha$ then $\sup_{t\le
N}P^{\eta_k}_{t}(y,\Gamma)< \epsilon $.}
$$ So if we take the limit (of the adequate subsequence) then the fact required
 is obtained.

For the second condition, we recall that throughout the proof of
Lemma 3.3 we obtained the following inequality :
$$\text{  $\forall t>0$, $\forall \eta_k>0$,
$\sup_{p\in K_D}\frac {Prob\{Y_t^{\eta_k}\in B_D^c(p,\epsilon)\}}{
\frac{t}{\eta_k^2}} \le
\frac{e^{-\frac{t}{\eta_k^2}}}{\frac{t}{\eta_k^2}}$.}$$

Then it is enough to let $\eta_k$ and at the same time $t$ go to
zero, to obtain the second condition, which ends the proof. \hfill
$\square$
\enddemo

\head{ 5.  Recurrent and transient behavior of the Brownian
motion}
\endhead

Usually, in the literature about the Brownian motion in the smooth
case, the authors question the {\it recurrent} or {\it transient}
behavior of this stochastic process.

It is known, for example, that the euclidian Brownian motion is
recurrent when it is two dimensional and it is transient if its
dimension is greater or equal to three. Moreover, we know that the
noncompact hyperbolic surface valued Brownian motion is transient.

For more results and details, we recommend to the reader the
papers of {\it H.P Mckean, D. Sullivan} [23] and {\it T.J. Lyons,
H.P Mckean} [21].

\vskip .5cm

\subhead{ 5.1 The geometric behavior of the admissible Riemannian
complex valued Brownian motion }\endsubhead

Let $K$ denotes a complete admissible Riemannian complex of
dimension $n$ and $p\in K$. We recall that the $K$-valued Brownian
motion
 $\{B^p_t\}_{t\ge0}$, was obtained as a weak limit of sequence
  of isotropic transport processes.

 On the other hand, we have seen that the trajectories of the
 isotropic processes are concatenations of geodesic segments. When
 a trajectory joins $(a.e.)$ transversally the
 $(n-1)$-skeleton$\setminus(n-2)$-skeleton, it goes on choosing isotropically a new maximal face
  (i.e. all adjacent maximal faces have the same probability to be chosen).

 Consequently, the $K$-valued Brownian motion $\{B^p_t\}_{t\ge0}$ behaves,
 inside every $n$-simplex $\Delta_{n}$, as the standard Brownian
motion with values in Riemannian $n$-dimensional manifold endowed
with the metric $g_{\Delta_{n}}$.

 Moreover, the process hits $(a.e.)$ "transversally" the
 $(n-1)$-skeleton$\setminus(n-2)$-skeleton, then it goes on choosing
 isotropically a maximal face.
 Thus, it results from this geometric description a new discreet random walk
 corresponding to the isotropic choices of the maximal faces.

To give a rigorous mathematical construction of such discreet
process, let us consider some notions.

 The {\it dual graph } $X$ of a complex
 $K$ is $1$-dimensional simplicial complex defined as follows:

 Consider one point inside (topological interior) each $n$-simplex of $K$ and, for every
  $(n-1)$-simplex a point in its topological interior, then, we connect the considered points with
  geodesic segments and let $E(X)$ denote the set of such segments. Thus, this dual graph consists
  of set $V_n(X)$ of vertices of degree $n+1$ (interior $n$-simplexes points) and a set $V_{n-1}(X)$
  of vertices corresponding to the interior points of the $(n-1)$-simplexes,
  where every vertex has degree equal to the number of the $n$-simplexes adjacent to this vertex.

  Consider now, the Markov chain (discreet Markov process) $\{C_n\}_{n\in \Bbb
  N}$ which has as a transition probability the function :
$$p(x,y) = \cases \frac{1}{deg\, x}&\text{ if $x,y \in V_{n-1}(X)$
and there exists $z\in V_n(X)$ such that $xz$, $zy \in E(X)$,}\\
 0 &\text{ unless,}\endcases$$ where $deg\, x$ is the degree of $x$ and $xz$ is an edge (geodesic segment)
 connecting $x$ to $z$.

Thus the latest random walk is a discreet "jump" process on the
set $V_{n-1}(X)$.

\vskip .5cm

\subhead{ 5.2 Brownian motion in an admissible complex with
nonpositive curvature and with dimension at the most
$2$}\endsubhead

This subsection is devoted to the study of the transient or
recurrent behavior of the Brownian motion in an admissible complex
with nonpositive curvature (in the sense of Alexandrov) and with
dimension at the most $2$.

Now recall the definition of recurrent/transient process:

\definition{Definition 5.2}

Let $\{X_t^p\}_t$ denote a stochastic process in a metric space
 $K$. Then $\{X_t^p\}_t$ is said to be recurrent if, for every ball
$B_p$ containing the point $p$, the process $\{X_t^p\}_t$ returns
to the ball $B_p$ (and so infinitely) with probability equal to
one ; in other words, the process is transient.

\enddefinition

\proclaim{Remark 5.3}

When the space $K$ is a discreet space, we consider the point $p$
instead of the ball $B_p$ in the above definition.
\endproclaim

\proclaim{Theorem 5.4}

Let $K$ denote a $2$-dimensional (respectively $1$-dimensional)
non-compact complete simply connected  admissible Riemannian
complex with nonpositive curvature. Then, if for every $1$-simplex
(respectively a vertex) there is at least three $2$-simplices
(respectively $1$-simplices) adjacent to it, the Brownian motion
is transient.
\endproclaim

Before proceeding to look at the proof of Theorem 5.4, let us
first give a short treatise on {\it simple random walk} on a
graph.

 Let $X= (V(X),E(X))$ denote a connected locally finite graph
 (a $1$-dimensional admissible Riemannian complex),
 where $V(X)$ is the set of vertexes and $E(X)$ is the set of
 edges.
 By {\it simple random walk} on the graph $X$,
 we mean the Markov chain for which the transition probability $p(x,y)$
 from vertex $x$ to vertex $y$ is given by the function :
 $$p(x,y) = \cases \frac{1}{deg\, x}, &\text{ if $xy \in E(X)$,}\\
 0, &\text{ unless},\endcases$$where $xy$ is an edge connecting
 $x$ to $y$.

 We say that $X$ is recurrent (respectively transient)
 if the simple random walk is recurrent (respectively
 transient).

 The {\it word metric} on the graph $X$ is an intrinsic metric in which
 each edge has unit length.

\proclaim{Remark 5.5 [13, Ch 6]}

Let $X$ denote a connected locally finite graph with uncountably
many ends. Assuming that every vertex has degree greater or equal
to three, then $X$ is transient.

\endproclaim

 \demo{Proof of Theorem 5.4}

 Let $K$ be an admissible complex  and let $X$ denote the dual graph of $K$. Now in the following,
 we will construct a graph $Y$ from the graph $X$.

 Let $x_1$ be a vertex belonging  to the set $V_{1}(X)$ and
 $z_1\in V_2(X)$ such that $x_1z_1\in E(X)$. Recall that the
 degree of $z_1$ is equal to three. We delete an edge adjacent to $z_1$,
 different than $x_1z_1$. We do the same thing with the other faces adjacent to
 $x_1$.

 Now go back to $z_1$, it is connected to another vertex
 $x_2 \in V_1(X)$ ($x_1$, $z_1$ and $x_2$ are all in the same $2$-simplex).
 We do the same thing with $x_2$ as we have done with $x_1$.
 At the end of this construction, forgetting the vertexes of
 degree equal to two, and as a consequence of the hypothesis on the complex $K$,
  we get a graph $Y$ isometrically equivalent to connected locally finite graph with uncountably
many ends and whose each vertex degree is greater or equal to
three. Moreover, the random walk coming from the isotropic choice
of maximal faces by the Brownian motion induces a simple random
walk on the graph $Y$.

Just now, suppose that the $K$-valued Brownian motion
$\{B^p_t\}_{t\ge 0}$ is recurrent. We can suppose that the point
$p$ is in the interior of an edge. Take as compact neighborhood of
the point $p$ the union of all its adjacent $2$-simplices and note
this neighborhood $B_p$.

Thus, if $\{B^p_t\}_{t\ge 0}$ returns to the ball $B_p$  with
probability equal to one, then inevitably, the simple random walk
on $Y$ returns to the point $p$ with probability $1$. In other
words, the graph $Y$ is recurrent which contradicts Remark 5.5,
and so the theorem is now proven. \hfill $\square$
\enddemo


\Refs \widestnumber\key{HPS}


     \ref \key 1 \by S.B.  Alexander, R.L.  Bishop \paper The Hadamard-Cartan
       theorem in locally convex metric spaces \paperinfo L'Enseignement Math, 36
                             , 309-320, (1990) \endref

    \ref \key 2 \by A.D.  Alexandrov \paper A theorem on triangles in a metric
 space and some applications \paperinfo Trudy Math.  Inst.Steklov 38,
 5-23, (Russian) (1951) \endref


 \ref \key 3 \by W.  Ballmann, M. Brin \paper Orhihedra of Nonpositive Curvature
                    \paperinfo Publications IHES , 82, 169-209, (1995)\endref

 \ref \key 4 \by W.  Ballmann, S. Buyalo \paper  Nonpositively Curved Metrics on 2-Polyhedra
                    \paperinfo Math. Zeitschrift, 222, No.1, 97-134 (1996)\endref
\ref \key 5 \by P. Billingsley \book Convergence of probability
Mesures     \bookinfo Wiley Series in Probability and Mathematical
Statistics , (1968)
\endref

 \ref \key 6 \by M.R. Bridson \paper Geodesics and Curvature in Metric Simplicial Complexes \paperinfo
     World Scientific, Eds. E. Ghys, A.Haefliger, A. Verjovsky, (1990)\endref

\ref\key 7 \by M.R. Bridson, A. Haefliger \book  Metric spaces of
Non-positive curvature \bookinfo Springer (1999) \endref





\ref \key 8 \by M.  Brin, Y.  Kifer \paper Brownian motion,
harmonic functions and hyperbolicity for Euclidean complexes
\paperinfo Math. Zeitschrift, 237, 421-468, (2001)
\endref

 \ref \key 9 \by H.  Busemann\paper Spaces with nonpositive curvature
               \paperinfo Acta Math. 80, 259-310, (1948) \endref

 \ref \key 10 \by G. De Cecco, G. Palmieri \paper Distanza
 intrinseca una variet\`a finsleriana di Lipschitz
               \paperinfo Rend. Aca. Naz. Sci. 17, 129-151, (1993) \endref

\ref \key 11 \by  M. Davis, T. Januzkiewicz\paper Hyperbolization
of polyhedra\paperinfo Journal of Differential Geometry, 34(2),
347-388, (1991) \endref



\ref \key 12 \by J.L. Doob \book Stochastic Processes
                    \bookinfo Wiley, New York, (1953) \endref

\ref \key 13 \by P.G. Doyle, J.L. Snell \book Random walks and
electric networks  \bookinfo Carus Math. Monographs, 22, Math.
Assoc. Amer., Washington, DC, (1984)
\endref



 \ref \key 14 \by E.B. Dynkin \book Markov Processes I, II
                    \bookinfo Springer-Verlag, (1965) \endref

\ref \key 15 \by J. Eells, B. Fuglede \book Harmonic maps between
Riemannian polyhedra \bookinfo  Cambridge university press, (2001)
\endref

\ref \key 16 \by M.  Emery \book Stochastic Calculus in Manifolds
 \bookinfo  Springer-Verlag Berlin Heidelberg New York, (1980) \endref

   \ref \key 17 \by E.  Ghys, P.  de la Harpe (ed) \paper Sur les groupes
      hyperboliques d'apr\`es M.  Gromov  \paperinfo Progress in Math.  83,
                              Birkhauser(1990) \endref
\ref\key 18\by M. Gromov \paper Structures m\'etrique pour les
vari\'et\'es Riemanniennes, r\'edig\'e par J.Lafontaine et P.Pansu
\paperinfo Cedic/Fernand, Nathan (1981) \endref



     \ref \key 19 \by E. J\o rgensen \paper The Central Limit Problem for Geodesic Random Walks \paperinfo
     Z Wahrscheinlichkeistheorie verw. Gebiete 32, 1-64, (1975) \endref

\ref \key 20 \by T. G. Kurtz \paper A Limit theorem for perturbed
operator semigroups with applications to random evolution
\paperinfo J. Functional Analysis 12, 55-67, (1973)
\endref


\ref \key 21 \by T. J. Lyons, H. P. Mckean \paper Winding of plane
Brownian motion \paperinfo
     Adv. in Math. 51, 212-225, (1984)\endref

\ref \key 22 \by P. Malliavin \paper Diffusions et g\'eom\'etrie
diff\'erentielle globale \paperinfo
     Lecture Notes, August, (1975)\endref


\ref \key 23 \by H. P. Mckean, D. Sullivan \paper Brownian motion
and Harmonic functions on the class surface of thrice punctured
sphere \paperinfo Adv. in Math. 51, 203-211, (1984)
    \endref



    \ref \key 24 \by J.  Milnor \book Morse Theory\bookinfo Princeton University Press
(1969) \endref

 \ref \key 25 \by J.  Munkres \book  TOPOLOGY a first course
\bookinfo  PENTICE HALL, Engelwood Cliffs, New Jersey (1975)
\endref


\ref \key 26 \by F. Paulin \paper Constructions of hyperbolic
groups via hyperbolizations of polyhedra \paperinfo  World
Scientific, Eds. E. Ghys,
 A. Haefliger, A. Verjovsky, (1990)\endref


   \ref \key 27 \by M.  Pinsky \paper Isotropic Transport Process On Riemannian Manifold
 \paperinfo Transactions Of American Mathematical Society 218, 353-360, (1976) \endref




\ref \key 28 \by A. V. Skorohod \paper Limit Theorem for
Stochastic processes  \paperinfo Teor. Verojatnost. i Primenen. 1,
261-290, (1956)
\endref

\ref \key 29 \by C. Stone \paper Weak Convergence of Stochastic
processes Defined On Semi-Infinite Time Intervals \paperinfo Proc.
Amer. Math. Soc. 14, 694-696, (1963)
\endref

 \ref \key 30 \by J.  Tits \book  Buildings of
spherical type and finite BN-pairs \bookinfo  volume 386 Springer,
 (1974)\endref


 \ref \key 31 \by S. Watanabe, T. Watanabe \paper  Convergence of isotropic scattering transport
 process to Brownian motion \paperinfo  Nagoya Math. J. 40, 161-171, (1970)\endref



                                      \endRefs
                                       \enddocument
\end